# SOME REMARKS ON FIRST PASSAGE OF LEVY PROCESSES, THE AMERICAN PUT AND PASTING PRINCIPLES

BY L. ALILI AND A. E. KYPRIANOU

*University of Warwick and Utrecht University*

The purpose of this article is to provide, with the help of a fluctuation identity, a generic link between a number of known identities for the first passage time and overshoot above/below a fixed level of a Lévy process and the solution of Gerber and Shiu [*Astin Bull.* **24** (1994) 195–220], Boyarchenko and Levendorskiĭ [Working paper series EERS 98/02 (1998), Unpublished manuscript (1999), *SIAM J. Control Optim.* **40** (2002) 1663–1696], Chan [Original unpublished manuscript (2000)], Avram, Chan and Usabel [*Stochastic Process. Appl.* **100** (2002) 75–107], Mordecki [*Finance Stoch.* **6** (2002) 473–493], Asmussen, Avram and Pistorius [*Stochastic Process. Appl.* **109** (2004) 79–111] and Chesney and Jeanblanc [*Appl. Math. Fin.* **11** (2004) 207–225] to the American perpetual put optimal stopping problem. Furthermore, we make folklore precise and give necessary and sufficient conditions for smooth pasting to occur in the considered problem.

**1. Introduction.** Let $X = \{X_t : t \geq 0\}$ be a Lévy process defined on filtered probability space $(\Omega, \mathcal{F}, \{\mathcal{F}_t\}, \mathbb{P})$ satisfying the usual conditions. For $x \in \mathbb{R}$, denote by $\mathbb{P}_x(\cdot)$ the law of $X$ when it is started at $x$ and write simply $\mathbb{P}_0 = \mathbb{P}$. We denote its Lévy–Khintchine exponent by $\Psi$, that is, $\mathbb{E}[e^{i\theta X_1}] = \exp\{-\Psi(\theta)\}$, $\theta \in \mathbb{R}$. Now, consider the following optimal stopping problem:

$$v(x) = \sup_{\tau \in \mathcal{T}_{0,\infty}} \mathbb{E}_x[e^{-r\tau}(K - e^{X_\tau})^+], \tag{1}$$

where $K > 0, r \geq 0$ and $\mathcal{T}_{0,\infty}$ is the family of optimal stopping times with respect to $\{\mathcal{F}_t : t \geq 0\}$. Establishing the optimal value and optimal stopping time in (1) is closely related to finding the value pricing and exercise strategy of an American put option in an incomplete Black–Scholes-type markets driven by Lévy processes (see [13, 35] or [18]). For this reason, we refer to (1) as the *American put optimal stopping problem*.

In a number of numerical simulations and theoretical calculations for specific choices of Lévy processes, various authors have found that the American put optimal stopping problem is solved in the same way as for the case that $X$ is a scaled Brownian motion with drift (the Black–Scholes market). In other words, it is solved by a strategy of the form

$$\tau^* = \inf\{t \geq 0 : X_t < x^*\}$$

for a specific value $x^* < \log K$ so that

$$v(x) = K\mathbb{E}_x[e^{-r\tau^*}] - \mathbb{E}_x[e^{-r\tau^* + X_{\tau^*}}],$$

thus, linking the American perpetual put optimal stopping problem to the first passage problem of a Lévy process. See Gerber and Shiu [21], Chan [15, 16], Boyarchenko and Levendorskiĭ [10, 11, 12, 13, 14], Mordecki [27, 28], Avram, Chan and Usabel [5], Asmussen, Avram and Pistorius [4], Chesney and Jeanblanc [17], Hirsa and Madan [23], Matache, Nitsche and Schwab [26], Almendral and Oosterlee [2] and Almendral [1]. Notably, Mordecki [28] handles the case when $X$ is a general Lévy process.

In addition, these authors have observed that the function $v$ is continuous, equals $K - e^x$ for $x \leq x^*$ and is bounded below by $(K - e^x)^+$ for $x > x^*$. Recall that $v$ is said to exhibit smooth pasting when $v$ pastes onto the lower bounding curve $(K - e^x)^+$ at $x^*$ in such a way that the left and right derivative there agree. Among the aforementioned list of articles, some authors have also observed that, unlike the case when $X$ is purely Gaussian, the optimal value function $v$ does not necessarily respect the so-called *principle of smooth pasting*. Boyarchenko and Levendorskiĭ [10, 11, 12, 13] supplied sufficient conditions for a class of Lévy processes with exponential moments and "stable" like characteristic exponent (RLPE processes) under which continuous or smooth pasting appears. For their class of process, they show that there is smooth pasting if $X$ has bounded variation with a nonpositive drift term or if $X$ has unbounded variation. Their method also offers identities for the value function and optimal exercise barrier for a more general claim structure and for Lévy processes which possess absolutely continuous resolvents; see [12, 14]. Chan [15, 16] also observed for the case of spectrally negative processes that there was no smooth pasting if and only if the process is of bounded variation.

In related work Peskir and Shiryaev [32, 33] (see also [31]) and Gapeev [20] studied a number of optimal stopping problems for special classes of Markov process of bounded variation with jumps such that the inter-arrival times of the jumps are independent and exponentially distributed. The optimal strategies in these problems consist of stopping when the underlying Markov process crosses specified thresholds. Further, one sees in these papers that the principle of smooth pasting could also not be taken for granted at the optimal stopping thresholds. However, as in the case of the American put option, continuity was still observed when there was no smooth pasting. In their papers, Peskir and Shiryaev give an intuitive explanation as to precisely when and why smooth pasting does not hold for their models and argue for a *principle of continuous pasting* instead.

The aim of this paper is to comprehensively link a variety of identities for first passage problems of different Lévy processes and their connection to the American put optimal stopping problem which have appeared in recent literature and to explain precisely when smooth pasting is absent. We shall do this with the help of a fluctuation identity. Mordecki [28] gave a closed form for the solution to (1) for a general Lévy process in terms of the distribution of the infimum of $X$ at an independent exponential time. The fluctuation identity will also allow us to give another proof of this form of the solution to (1), avoiding the need for a random walk approximation (as was the case in Mordecki's original proof) and giving insight into the importance of the role played by the regularity of the paths of $X$ in the solution. Further, we give necessary and sufficient conditions for the solution to exhibit smooth pasting. It turns out that the principle of smooth pasting boils down to, quite simply, the regularity of the point 0 for $(-\infty, 0)$ for $X$.

The paper closes with a conjecture. We believe this can, in principle, be used as a rule of thumb to predict for Markovian optimal stopping problems whether the principle of smooth pasting holds.

**2. A fluctuation identity for overshoots.** An important tool in the study of the fluctuations of Lévy processes is the Wiener–Hopf factorization which we now review for convenience. For a more detailed account, the reader is referred to [22] or [7].

Assume that $r > 0$. We have that $\overline{X}_{\mathbf{e}_r}$ and $\overline{X}_{\mathbf{e}_r} - X_{\mathbf{e}_r}$ are independent where $\overline{X}_t = \sup_{0 \le s \le t} X_t$ and $\mathbf{e}_r$ is an independent exponentially distributed random variable with parameter $r$. As a consequence, for $\theta \in \mathbb{R}$,

$$\frac{r}{r + \Psi(\theta)} = \Psi_r^+(i\theta) \cdot \Psi_r^-(i\theta)$$

with

(2) $$\Psi_r^+(i\theta) = \mathbb{E}[e^{i\theta \overline{X}_{\mathbf{e}_r}}]$$

and

(3) $$\Psi_r^-(i\theta) = \mathbb{E}[e^{-i\theta(\overline{X}_{\mathbf{e}_r} - X_{\mathbf{e}_r})}] = \mathbb{E}[e^{i\theta \underline{X}_{\mathbf{e}_r}}],$$

where $\underline{X}_t = \inf_{0 \le s \le t} X_s$. The second equality in (3) follows by duality, that is, when $X$ is time reversed over a fixed time interval, it has the same law as $-X$. Note that $\Psi^+$ can be analytically extended to $\{z \in \mathbb{C} : \Re z \leq 0\}$ and $\Psi^-$ can be analytically extended to $\{z \in \mathbb{C} : \Re z \geq 0\}$.

We present a fluctuation identity which looks, at a first sight, trivial. However, with the help of the Wiener–Hopf factorization, the Pecherskii–Rogozin identity follows from it in a very straightforward way. As indicated in the Introduction, the identity also gives clarification to a number of other existing identities and their role played in the solution to the optimal stopping problem (1). Essentially, this identity is not new, as it appears implicitly in a number of texts dating back to at least the seventies. See, for example, [19], page 1368, where one sees the same identity for random walks embedded in the proof of another result.

To a fixed level $x \in \mathbb{R}$ we associate the first strict passage time $\tau_x^+$ (resp. $\tau_x^-$) above (resp. below) $x$, that is,

$$\tau_x^+ = \inf\{t \geq 0 : X_t > x\} \quad \text{and} \quad \tau_x^- = \inf\{t \geq 0 : X_t < x\}.$$

Now, we are ready to state the identity.

LEMMA 1. *For all $\alpha > 0$, $\beta \geq 0$ and $x \geq 0$, we have*

$$\mathbb{E}[e^{-\alpha \tau_x^+ - \beta X_{\tau_x^+}} \mathbf{1}_{(\tau_x^+ < \infty)}] = \frac{\mathbb{E}[e^{-\beta \overline{X}_{\mathbf{e}_\alpha}} \mathbf{1}_{(\overline{X}_{\mathbf{e}_\alpha} > x)}]}{\mathbb{E}[e^{-\beta \overline{X}_{\mathbf{e}_\alpha}}]}. \tag{4}$$

*For the case that $\alpha = 0$, the left-hand side has a limit and, hence, so does the right-hand side. If we additionally assume that $\overline{X}_\infty < \infty$, then we may identify the limit as follows:*

$$\mathbb{E}[e^{-\beta X_{\tau_x^+}} \mathbf{1}_{(\tau_x^+ < \infty)}] = \frac{\mathbb{E}[e^{-\beta \overline{X}_\infty} \mathbf{1}_{(\overline{X}_\infty > x)}]}{\mathbb{E}[e^{-\beta \overline{X}_\infty}]}.$$

PROOF. First, assume that $\alpha, \beta, x > 0$ and note that

$$\mathbb{E}[e^{-\beta \overline{X}_{\mathbf{e}_\alpha}} \mathbf{1}_{(\overline{X}_{\mathbf{e}_\alpha} > x)}] = \mathbb{E}[e^{-\beta \overline{X}_{\mathbf{e}_\alpha}} \mathbf{1}_{(\tau_x^+ < \mathbf{e}_\alpha)}]$$

$$= \mathbb{E}[\mathbf{1}_{(\tau_x^+ < \mathbf{e}_\alpha)} e^{-\beta X_{\tau_x^+}} \mathbb{E}[e^{-\beta(\overline{X}_{\mathbf{e}_\alpha} - X_{\tau_x^+})} | \mathcal{F}_{\tau_x^+}]].$$

Now, conditionally on $\mathcal{F}_{\tau_x^+}$ and on the event $\tau_x^+ < \mathbf{e}_\alpha$, the random variables $\overline{X}_{\mathbf{e}_\alpha} - X_{\tau_x^+}$ and $\overline{X}_{\mathbf{e}_\alpha}$ have the same distribution, thanks to the lack of memory property of $\mathbf{e}_\alpha$ and the strong Markov property. Hence, we have the factorization

$$\mathbb{E}[e^{-\beta \overline{X}_{\mathbf{e}_\alpha}} \mathbf{1}_{(\overline{X}_{\mathbf{e}_\alpha} > x)}] = \mathbb{E}[e^{-\alpha \tau_x^+ - \beta X_{\tau_x^+}}] \mathbb{E}[e^{-\beta \overline{X}_{\mathbf{e}_\alpha}}].$$

The case that $\alpha$, $\beta$ or $x$ are equal to zero can be achieved by taking limits on both sides of the above equality. □

By replacing $X$ by $-X$ in the previous lemma, we get the following result.

COROLLARY 2. *For all $\alpha, \beta \geq 0$ and $x \geq 0$, we have*

$$\tag{5} \mathbb{E}[e^{-\alpha \tau_{-x}^- + \beta X_{\tau_{-x}^-}} \mathbf{1}_{(\tau_{-x}^- < \infty)}] = \frac{\mathbb{E}[e^{\beta \underline{X}_{\mathbf{e}_\alpha}} \mathbf{1}_{(-\underline{X}_{\mathbf{e}_\alpha} > x)}]}{\mathbb{E}[e^{\beta \underline{X}_{\mathbf{e}_\alpha}}]}.$$

**3. First passage.** In this section we shall use (1), together with the Wiener–Hopf factorization, to recover the Pecherskii–Rogozin identity (sometimes called the second factorization identity). Further, we shall cross reference (4) and (5) against the collection of explicitly or semi-explicitly known first passage identities.

3.1. *The Pecherskii–Rogozin identity.* By taking Laplace transforms of both sides of (4) and using Fubini's theorem, we can write, for $q > 0$,

$$\int_0^\infty e^{-qx} \mathbb{E}[e^{-\alpha \tau_x^+ - \beta(X_{\tau_x^+} - x)}] \, dx$$

$$= \frac{1}{\mathbb{E}[e^{-\beta \overline{X}_{\mathbf{e}_r}}]} \int_0^\infty e^{-qx} \mathbb{E}[e^{-\beta(\overline{X}_{\mathbf{e}_r} - x)} \mathbf{1}_{(\overline{X}_{\mathbf{e}_\alpha} > x)}] \, dx$$

$$= \frac{1}{\mathbb{E}[e^{-\beta \overline{X}_{\mathbf{e}_\alpha}}]} \int_0^\infty e^{-qx} \int_0^\infty \mathbf{1}_{(y > x)} e^{-\beta(y - x)} \mathbb{P}(\overline{X}_{\mathbf{e}_\alpha} \in dy) \, dx$$

$$= \frac{1}{\mathbb{E}[e^{-\beta \overline{X}_{\mathbf{e}_\alpha}}]} \int_0^\infty e^{-\beta y} \int_0^\infty \mathbf{1}_{(y > x)} e^{-qx + \beta x} \, dx \, \mathbb{P}(\overline{X}_{\mathbf{e}_\alpha} \in dy)$$

$$= \frac{1}{(q - \beta) \mathbb{E}[e^{-\beta \overline{X}_{\mathbf{e}_\alpha}}]} \int_0^\infty (e^{-\beta y} - e^{-qy}) \mathbb{P}(\overline{X}_{\mathbf{e}_\alpha} \in dy)$$

$$= \frac{\mathbb{E}[e^{-\beta \overline{X}_{\mathbf{e}_\alpha}}] - \mathbb{E}[e^{-q \overline{X}_{\mathbf{e}_\alpha}}]}{(q - \beta) \mathbb{E}[e^{-\beta \overline{X}_{\mathbf{e}_\alpha}}]}.$$

If $X$ is an ascending subordinator with Laplace exponent $\phi(q) := -\log \mathbb{E}(e^{-qX_1})$, $q \geq 0$, then we see that

$$\int_0^\infty e^{-qx} \mathbb{E}[e^{-\alpha \tau_x^+ - \beta(X_{\tau_x^+} - x)}] \, dx = \frac{\phi(q) - \phi(\beta)}{(q - \beta)(\alpha + \phi(q))},$$

which was also established in [37]. Otherwise, if $X$ is not a descending subordinator, with the help of (2), we come to rest at the Pecherskii–Rogozin identity

$$\int_0^\infty e^{-qx} \mathbb{E}[e^{-\alpha \tau_x^+ - \beta(X_{\tau_x^+} - x)}] \, dx = \frac{1}{q - \beta} \left(1 - \frac{\Psi_\alpha^+(-q)}{\Psi_\alpha^+(-\beta)}\right)$$

for any $q > 0$. See [30] and [34] for a comparison with existing proofs.

3.2. *Spectrally one-sided processes.* Suppose that $X$ is spectrally negative, but not a negative subordinator, with Laplace exponent $\psi\colon [0,\infty) \to \mathbb{R}$ given by $\mathbb{E}[\exp\{\lambda X_t\}] = \exp\{\psi(\lambda) t\}$ so that $\psi(z) = -\Psi(-iz)$ for $\Re z \geq 0$. Denote by $\Phi(\alpha)$ the largest root of the equation $\psi(\lambda) - \alpha = 0$, where $\alpha \geq 0$. The Wiener–Hopf factors are now given by

$$\Psi_\alpha^+(i\theta) = \frac{\Phi(\alpha)}{\Phi(\alpha) - i\theta} \quad \text{and} \quad \Psi_\alpha^-(i\theta) = \frac{\alpha}{\Phi(\alpha)} \frac{\Phi(\alpha) - i\theta}{\alpha - \psi(i\theta)}.$$

When working with Laplace transforms in (2) and (3), we have that

$$\int_{[0,\infty)} e^{-\beta x} \mathbb{P}(\overline{X}_{\mathbf{e}_\alpha} \in dx) = \frac{\Phi(\alpha)}{\Phi(\alpha) + \beta}$$

and

$$\int_{[0,\infty)} e^{-\beta x} \mathbb{P}(-\underline{X}_{\mathbf{e}_\alpha} \in dx) = \frac{\alpha}{\Phi(\alpha)} \frac{\beta - \Phi(\alpha)}{\psi(\beta) - \alpha}.$$

Therefore, we see the known distributional identities (cf. [7, 9]),

(6) $$\mathbb{P}(\overline{X}_{\mathbf{e}_\alpha} \in dx) = \Phi(\alpha) e^{-\Phi(\alpha) x}\, dx$$

and

(7) $$\mathbb{P}(-\underline{X}_{\mathbf{e}_\alpha} \in dx) = \frac{\alpha}{\Phi(\alpha)}\, dW^{(\alpha)}(x) - \alpha W^{(\alpha)}(x)\, dx,$$

where $W^{(\alpha)}\colon [0,\infty) \to [0,\infty)$ is the scale function (cf. [9]) and is characterized on $(0,\infty)$ by

$$\int_0^\infty e^{-\lambda x} W^{(\alpha)}(x)\, dx = \frac{1}{\psi(\lambda) - \alpha} \qquad \text{for } \lambda > \Phi(\alpha).$$

Plugging (6) and (7) into (4) and (5), respectively, we arrive at the well-known expression

$$\mathbb{E}[e^{-\alpha \tau_x^+ - \beta X_{\tau_x^+}} \mathbf{1}_{(\tau_x^+ < \infty)}] = e^{-(\Phi(\alpha) + \beta) x}.$$

After some algebra, we also get the accompanying one

$$\mathbb{E}[e^{-\alpha \tau_{-x}^- + \beta X_{\tau_{-x}^-}} \mathbf{1}_{(\tau_{-x}^- < \infty)}] = (\psi(\beta) - \alpha) \int_x^\infty e^{-\beta y} W^{(\alpha)}(y)\, dy$$

$$- \frac{\alpha - \psi(\beta)}{\Phi(\alpha) - \beta} e^{-\beta x} W^{(\alpha)}(x)$$

for all $\alpha, \beta, x \geq 0$.

3.3. *Subordinators.* Take $X$ to be a subordinator and denote by $\phi$ its Laplace–Lévy exponent satisfying $\mathbb{E}[e^{-\lambda X_t}] = e^{-t\phi(\lambda)}$ for $\lambda \geq 0$. We also introduce the $\alpha$-resolvent $U^{(\alpha)}$ defined via

$$\frac{1}{\alpha + \phi(\lambda)} = \int_0^\infty e^{-\alpha t}\mathbb{E}[e^{-\lambda X_t}]\,dt = \int_0^\infty e^{-\lambda y}U^{(\alpha)}(dy).$$

Because $\overline{X}_{\mathbf{e}_\alpha} = X_{\mathbf{e}_\alpha}$, identity (4) simply reads

(8) $$\mathbb{E}[e^{-\alpha \tau_x - \beta X_{\tau_x^+}}] = (\alpha + \phi(\beta))\int_{[x,\infty)} e^{-\beta z}U^{(\alpha)}(dz)$$

for all $\alpha, \beta, x \geq 0$.

3.4. *A general phase-type process.* Here we borrow an example which appeared in [29] and then with a different proof in [4]. The Lévy process $X$ is taken to be the independent sum of a spectrally positive process with a compound Poisson process having negative phase-type jumps.

Recall that a distribution $F$ on $(0, \infty)$ is *phase-type* if it is the distribution of the absorption time in a finite state continuous time Markov process $\{J_t : t \geq 0\}$ with one state $\Delta$ absorbing and the remaining ones $1, \ldots, m$ transient. The parameters of this system are $m$, the restriction $\mathbf{T}$ of the full intensity matrix to the $m$ transient states and the initial probability (row) vector $\mathbf{a} = (a_1, \ldots, a_m)$, where $a_i = \mathbb{P}(J_0 = i)$. For any $i = 1, \ldots, m$, let $t_i$ be the intensity of the transition $i \to \Delta$ and write $\mathbf{t}$ for the column vector of intensities. It follows that $F$ has a density given by $f(x) = \mathbf{a}e^{\mathbf{T}x}\mathbf{t}$ and Laplace transform given by $\widehat{F}(s) = \int_0^\infty e^{-sx}f(x)\,dx = \mathbf{a}(s\mathbf{I} - \mathbf{T})^{-1}\mathbf{t}$ for $s > 0$. The latter can be extended to the complex plane except at a finite number of poles (the eigenvalues of $\mathbf{T}$). For full details, we refer the reader to [3]. The process $X$ enjoys the representation

$$X_t = X_t^{(+)} - \sum_{j=1}^{N(t)} U_j, \qquad t \geq 0,$$

where $\{X_t^{(+)} : t \geq 0\}$ is a spectrally positive Lévy process, $\{N_t : t \geq 0\}$ is a Poisson process with rate $\lambda$ and $\{U_j : j \geq 1\}$ are i.i.d. random variables with a common distribution $F$; all of the aforementioned objects being mutually independent.

The corresponding Lévy–Khintchine exponent, $\Psi$, can be analytically extended to $\{z \in \mathbb{C} : \Re z \leq 0\}$ with the exception of a finite number of poles (the eigenvalues of $\mathbf{T}$). Define, for each $\alpha > 0$, the finite set of roots with negative real part

$$\mathcal{I}_\alpha = \{\rho_i : \Psi(\rho_i) + \alpha = 0, \Re\rho_i < 0\},$$

where multiple roots are counted individually. Next, define, for each $\alpha > 0$, a second set of roots with negative real part

$$\mathcal{J}_\alpha = \left\{ \eta_i : \frac{1}{\Psi(\eta_i) + \alpha} = 0, \Re \eta_i < 0 \right\},$$

again taking multiplicity into account. Mordecki [29] and Asmussen, Avram and Pistorius [4] show that

$$\Psi_\alpha^-(s) = \frac{\prod_{i \in \mathcal{I}_\alpha}(-\rho_i)}{\prod_{i \in \mathcal{J}_\alpha}(-\eta_i)} \frac{\prod_{i \in \mathcal{I}_\alpha}(s - \eta_i)}{\prod_{i \in \mathcal{J}_\alpha}(s - \rho_i)},$$

which can be analytically extended to the whole of $\mathbb{C}$ except for the poles at $\{\rho_j \in \mathcal{I}_\alpha\}$. Further, they show that

$$\mathbb{P}(-\underline{X}_{\mathbf{e}_\alpha} \in dx) = \sum_{j=1}^n \sum_{k=1}^{m(j)} A_{j,k} \frac{(-\rho_j x)^{k-1}}{(k-1)!} e^{\rho_j x} \, dx,$$

where $m(j)$ is the multiplicity of root $\rho_j$, $n$ is the number of distinct roots and

$$A_{j,k} = \frac{1}{(m-k)!} \frac{d^{m-k}}{ds^{m-k}} \left( \frac{\Psi_\alpha^-(s)(s - \rho_j)^m}{(-\rho_j)^k} \right) \bigg|_{s=\rho_j}.$$

According to (5), we now have

$$\mathbb{E}[e^{-\alpha \tau_{-x}^- + \beta X_{\tau_{-x}^-}} \mathbf{1}_{(\tau_{-x}^- < \infty)}]$$
$$= \sum_{j=1}^n \sum_{k=1}^{m(j)} A_{j,k} (-\rho_j)^{k-1} e^{(\rho_j - \beta)x} \sum_{i=0}^{k-1} \frac{1}{(k-1-i)!} \frac{x^{k-1-i}}{(\beta - \rho_j)^{i+1}}.$$

When none of the roots $\rho_j$ are multiple, we may simply write $A_j$ in place of $A_{j,1}$ and the previous formula is exactly the same as the one given in equation (24), Proposition 2 of [4].

**4. American put.** The following result was established for special, but nonetheless, rich classes of Lévy processes in [10, 11, 12, 13, 15, 16] and for a general Lévy process by Mordecki [28]. We give the version as it appears in [28].

THEOREM 3. *Assume that either $r > 0$ or $r = 0$ and*

$$\mathbb{P}\left(\lim_{t \uparrow \infty} X_t = \infty\right) = 1.$$

*Then*

$$v(x) = \frac{\mathbb{E}[(K\mathbb{E}[e^{\underline{X}_{\mathbf{e}_r}}] - e^{x + \underline{X}_{\mathbf{e}_r}})^+]}{\mathbb{E}[e^{\underline{X}_{\mathbf{e}_r}}]}$$

*and the optimal stopping time is given by*

$$\tau^* = \inf\{t \geq 0 : X_t < x^*\},$$

*where*

$$x^* = \log K \mathbb{E}[e^{\underline{X}_{\mathbf{e}_r}}].$$

In the sequel we shall work with the definition that $\mathbf{e}_r = \infty$ almost surely when $r = 0$.

Mordecki's proof consisted of first solving the analogue of (1) in the random walk setting. This was done in the spirit of reasoning found in [19] for a closely related optimal stopping problem for random walks based on the gain function $(e^x - K)^+$. Having obtained the result for random walks, Mordecki established the result for Lévy processes by checking that the result for random walks "converges" to the required identity for Lévy processes when passing to a limit.

REMARK 4. Note that the above theorem would appear to be incomplete. Since Lévy processes either drift to plus or minus infinity or oscillate, the case that $r = 0$ and $\liminf_{t \uparrow \infty} X_t = -\infty$ is missing. However, this case is not interesting for the current discussion, as it is easy to show that it is not optimal to stop in a finite time.

REMARK 5. Since for each positive $a, b$ the function $y \mapsto (a - by)^+$ is convex, it follows that the optimal value function $v$ is convex in the variable $y = e^x$ and, hence, $v$ is continuous and has left and right derivatives in $x$. See, for example, the discussion in [25].

In Section 6 we show how the identity given in Corollary 2 can be used to give a proof of Theorem 3 which avoids the necessity of a random walk approximation. Embedded in the proof is a clearer indication of the role played by the relevant pasting conditions at the optimal value $x^*$. In addition, the following theorem establishes necessary and sufficient conditions under which smooth pasting occurs.

THEOREM 6. *If either $r > 0$ or $r = 0$ and $\mathbb{P}(\lim_{t \uparrow \infty} X_t = \infty) = 1$, then the right derivative at the optimal stopping boundary is given by $v'(x^*+) = -e^{x^*} + K\mathbb{P}(\underline{X}_{\mathbf{e}_r} = 0)$. Thus, the optimal stopping problem (1) exhibits smooth pasting at $x^*$ if and only if $0$ is regular for $(-\infty, 0)$.*

The proofs of Theorems 3 and 6 are given in Section 6.

**5. Consistency with existing literature.**

5.1. *The optimal value function.* Using the examples given in Sections 3.2 and 3.4, together with the fact that Corollary 2 implies that

$$v(x) = K\mathbb{P}(-\underline{X}_{\mathbf{e}_r} \geq x - x^*) - Ke^{x-x^*}\mathbb{E}[e^{\underline{X}_{\mathbf{e}_r}}\mathbf{1}_{(-\underline{X}_{\mathbf{e}_r} \geq x-x^*)}]$$

(9)
$$= K\mathbb{E}[e^{-r\tau^-_{(x^*-x)}}] - e^x\mathbb{E}[e^{-r\tau^-_{(x^*-x)} + X_{\tau^-_{(x^*-x)}}}],$$

it is quite straightforward to verify that the expression for $v$ is consistent with that given in [4, 5, 15, 16, 21].

The class of regular Lévy processes of exponential type was introduced by Boyarchenko and Levendorskiĭ [12]. This class includes, for example, normal inverse Gaussian processes, hyperbolic processes and tempered stable processes, but not Variance Gamma processes. Boyarchenko and Levendorskiĭ [12] established an expression for the Fourier transform of the value function in (1) which involves one of the Wiener–Hopf factors. This is used to extract a sufficient condition for the case when there is no smooth pasting. We shall show below that, with straightforward calculations, one can derive the same expression for the Laplace transform of $v$, but not for a general Lévy process.

Using analytic extension, (2) and the Pecherskii–Rogozin identity, we can easily check that, for any $\lambda \in \mathbb{R}$,

(10)
$$\int_{x^*}^{\infty} e^{(i\lambda+\beta)(x-x^*)}\mathbb{E}[e^{\beta\underline{X}_{\mathbf{e}_r}}\mathbf{1}_{(-\underline{X}_{\mathbf{e}_r} > x-x^*)}]\,dx$$
$$= -\frac{1}{i\lambda+\beta}(\Psi_r^-(\beta) - \Psi_r^-(-i\lambda)).$$

Note when $r = 0$, this identity still makes sense, providing we keep to the assumption that $\lim_{t\uparrow\infty} X_t = \infty$. Now perform an elementary calculation to deduce that, for $\lambda \in \mathbb{R}$,

(11)
$$\int_{-\infty}^{x^*} e^{i\lambda x}(K - e^x)^+\,dx = \left(\frac{K}{i\lambda} - \frac{e^{x^*}}{i\lambda+1}\right)e^{i\lambda x^*}.$$

Observe that $v(x) = (K - e^x)^+$ for $x \leq x^*$, which combined with (9), (10) and (11), gives

$$\int_{\mathbb{R}} e^{i\lambda x} v(x)\,dx = -K\frac{1}{i\lambda}(1 - \Psi_r^-(-i\lambda))e^{i\lambda x^*}$$
$$+ K\frac{1}{i\lambda+1}(\Psi_r^-(1) - \Psi_r^-(-i\lambda))e^{i\lambda x^*}$$
$$+ \frac{K}{i\lambda}e^{i\lambda x^*} - \frac{e^{(i\lambda+1)x^*}}{i\lambda+1}.$$

Finally, noting that $e^{x^*} = K\Psi_r^-(1)$ in Theorem 3, we come to rest at the expression given in Section 4.2 of [12],

$$\int_{\mathbb{R}} e^{i\lambda x} v(x)\,dx = K \frac{e^{i\lambda x^*}}{i\lambda(i\lambda+1)} \Psi_r^-(-i\lambda)$$

for $\lambda \in \mathbb{R}$.

5.2. *Pasting condition.* To assist with the forthcoming analysis, let us recount some facts about the regularity of 0 for $(-\infty, 0)$ for Lévy processes. First, recall that 0 is regular for $(-\infty, 0)$ if and only if $\mathbb{P}(\tau_0^- = 0) = 1$. Thanks to Blumenthal's zero–one law, the latter probability equals 0 or 1. It is not difficult to construct an example of a Lévy process for which the $\mathbb{P}(\tau_0^- = 0) = 0$ (see Proposition 7 below). From the definition of regularity, one also readily sees that, under the assumptions of Theorem 3, when $r \geq 0$, $\underline{X}_{\mathbf{e}_r}$ has an atom at zero if and only if 0 is irregular for $(-\infty, 0)$; that is, to say $\mathbb{P}(\underline{X}_{\mathbf{e}_r} = 0) > 0$.

Let us now proceed to the description of the regularity of 0 for $(-\infty, 0)$ in terms of the underlying characteristics of $X$. The process $X$ has bounded variation if and only if

$$\int_{\mathbb{R}\setminus\{0\}} (1 \wedge |x|)\Pi(dx) < \infty,$$

where $\Pi$ is the characteristic Lévy measure, in which case we may write its Lévy–Khintchine exponent in the form

(12) $$\Psi(\theta) = -i\,\mathrm{d}\theta + \int_{\mathbb{R}\setminus\{0\}} (1 - e^{i\theta x})\Pi(dx).$$

The following proposition is a summary of what is known from the literature.

PROPOSITION 7. *The point 0 is regular for $(-\infty, 0)$ if and only if one of the following three conditions hold:*

  (i) *$X$ has bounded variation and $\mathrm{d} < 0$.*
 (ii) *$X$ has bounded variation, $\mathrm{d} = 0$ and*

$$\int_{-1}^{0-} \frac{|x|\Pi(dx)}{\int_0^{|x|} \Pi(y,\infty)\,dy} = \infty.$$

(iii) *$X$ has unbounded variation.*

Case (ii) was recently added to the class of processes exhibiting regularity of 0 for the lower half line in [8] and, for the other cases, we refer to the discussion at the beginning of [7], Section VI.3.

In [12] a sufficient condition for smooth pasting and a sufficient condition for continuous pasting alone were established for a class of Lévy processes possessing exponential moments and having Lévy–Khintchine exponent resembling that of a stable process. In terms of the notation given here, they show that there is smooth pasting if $X$ has bounded variation and $d \leq 0$ or if $X$ has unbounded variation. This suggests that the integral test in Proposition 7(ii) is automatically satisfied for bounded variation and $d = 0$ within the RLPE class. Taking, for example, the case of a CGMY- or KoBoL-type Lévy process, which is included in their class of Lévy processes, it is easy to confirm that the integral test is indeed automatically satisfied when $X$ has bounded variation and $d = 0$; see also the calculations in [24]. For a general RLPE Lévy process, however, to see why the integral test in Proposition 7 is satisfied, one may reason intuitively as follows. The small jump structure of an RLPE process is essentially the same as that of a symmetric stable process; this being due to a polynomial singularity in the Lévy measure close to zero. For the latter class of Lévy processes, the integral test, which concerns small jumps, is satisfied; hence, it is satisfied for the RLPE.

For the class of spectrally one sided models, Gerber and Shiu [21], who work with spectrally positive Lévy processes, derive an expression for $v$ from which one sees clearly that there is always smooth pasting. This is consistent with Theorem 6 since spectrally positive processes always have the required regularity property (cf. Corollary VII.5 of [7]). In [15, 16] it was shown for the case of spectrally negative Lévy processes (but not a negative subordinator) that there is smooth pasting if and only if $X$ has unbounded variation. This would seem to slightly differ from Theorem 6. However, we note that spectrally negative Lévy process of bounded variation (which are not negative subordinators) necessarily have a Lévy–Khintchine exponent in the form (12) with $d > 0$. Hence, according to Proposition 7, regularity for the lower half line coincides with unbounded variation. Note that similar conclusions for spectrally negative Lévy processes were also drawn for another optimal stopping problem related to the pricing of Russian options in [6].

Finally, within the class of compound Poisson models, Mordecki [27, 28] derives explicit formulae for the case of drifting Brownian motion, plus (mixed) exponential jumps which always has smooth pasting. Clearly, such processes always have unbounded variation on account of the presence of the Gaussian component. Hence, Theorem 6 exhibits consistency with Mordecki's findings.

**6. Proofs of Theorems 3 and 6.** The optimal stopping problem (1) has an infinite horizon and a Markovian claim structure. Based on many examples of optimal stopping problems with these characteristics, the proof of Theorem 3 works on the intuition that the optimal stopping problem is solved by a first passage time. Specifically, it appeals to the *classical* method of checking the following *well-known* sufficient conditions, which confirm that a *bounded* function $V$ and first passage time $\tau_y^-$ (for some $y \in \mathbb{R}$) characterize the value function and optimal stopping time of (1):

(i) $\{e^{-r(t \wedge \tau_y^-)} V(X_{t \wedge \tau_y^-}) : t \geq 0\}$ is a $\mathbb{P}_x$-martingale for all $x \in \mathbb{R}$,
(ii) $\{e^{-rt} V(X_t) : t \geq 0\}$ is a $\mathbb{P}_x$-supermartingale for all $x \in \mathbb{R}$,
(iii) $V(x) = (K - e^x)^+$ for all $x \leq y$,
(iv) $V(x) \geq (K - e^x)^+$ for all $x \in \mathbb{R}$.

To see why these are sufficient conditions, note that

$$V(x) \stackrel{\text{(ii)}}{\geq} \sup_{\tau \in \mathcal{T}_{0,\infty}} \mathbb{E}_x[e^{-r\tau} V(X_\tau)]$$
$$\stackrel{\text{(iv)}}{\geq} \sup_{\tau \in \mathcal{T}_{0,\infty}} \mathbb{E}_x[e^{-r\tau} (K - e^{X_\tau})^+]$$
$$\geq \mathbb{E}_x[e^{-r\tau_y^-} (K - e^{X_{\tau_y^-}})^+]$$
$$\stackrel{\text{(iii)}}{=} \mathbb{E}_x[e^{-r\tau_y^-} V(X_{\tau_y^-})]$$
$$\stackrel{\text{(i)}}{=} V(x),$$

showing that all inequalities are equalities. It follows that

$$V(x) = \sup_{\tau \in \mathcal{T}_{0,\infty}} \mathbb{E}_x[e^{-r\tau}(K - e^{X_\tau})^+] = \mathbb{E}_x[e^{-r\tau_y^-}(K - e^{X_{\tau_y^-}})^+].$$

The identity in Corollary 2 will play an instrumental role in finding the right pair $(V, \tau_y^-)$.

REMARK 8. To some extent, the conditions (i)–(iv) may be seen as a stochastic analogue of a free boundary value problem which one often sees as a way of characterizing the solution to an optimal stopping problem; see [36]. In this case it is also possible to write down a free boundary value problem (cf. [12]), although one must be a little careful about the sense in which the associated integro-differential operator is understood. Also, one must have a sense of the role played by the pasting condition in such a free boundary value problem; something we are currently investigating in this paper.

PROOF OF THEOREM 3. For each $y \in \mathbb{R}$, let us define the bounded functions $v_y(x) = \mathbb{E}_x[e^{-r\tau_y^-}(K - e^{X_{\tau_y^-}})^+]$. We shall reason our way to finding the right choice of $y$ for which $v_y$ satisfies (i)–(iv) above. By using Corollary 2,

we see immediately that

$$v_y(x) = \frac{\mathbb{E}[(K\mathbb{E}[e^{\underline{X}_{\mathbf{e}_r}}] - e^{x+\underline{X}_{\mathbf{e}_r}})\mathbf{1}_{(-\underline{X}_{\mathbf{e}_\alpha} > x - y)}]}{\mathbb{E}[e^{\underline{X}_{\mathbf{e}_r}}]}. \tag{13}$$

*Martingale property* (i). Begin by defining, for all $x \in \mathbb{R}$ and $\alpha, \beta \geq 0$, the function

$$h_{\alpha,\beta}(x) = \mathbb{E}_x[e^{-\alpha \tau_0^- + \beta X_{\tau_0^-}}] = \frac{e^{\beta x}\mathbb{E}[e^{\beta \underline{X}_{\mathbf{e}_\alpha}}\mathbf{1}_{(-\underline{X}_{\mathbf{e}_\alpha} > x)}]}{\mathbb{E}[e^{\beta \underline{X}_{\mathbf{e}_\alpha}}]},$$

where the second equality is a result of Corollary 2. Note, in particular, that, for $x \leq 0$, we have $h_{\alpha,\beta}(x) = \exp\{\beta x\}$, which, in combination with the strong Markov property, yields

$$\mathbb{E}_x[e^{-\alpha \tau_0^- + \beta X_{\tau_0^-}} | \mathcal{F}_t] = \mathbf{1}_{(t < \tau_0^-)} e^{-\alpha t} \mathbb{E}_{X_t}[e^{-\alpha \tau_0^- + \beta X_{\tau_0^-}}]$$
$$+ \mathbf{1}_{(t \geq \tau_0^-)} e^{-\alpha \tau_0^- + \beta X_{\tau_0^-}}$$
$$= e^{-\alpha(t \wedge \tau_0^-)} h_{\alpha,\beta}(X_{t \wedge \tau_0^-}).$$

Thus, $\{e^{-(t \wedge \tau_0^-)} h_{\alpha,\beta}(X_{t \wedge \tau_0^-}) : t \geq 0\}$ is a $\mathbb{P}_x$-martingale. Now note from (13) that $v_y$ is a linear combination of $h_{r,0}(x-y)$ and $h_{r,1}(x-y)$. It follows that $\{e^{-(t \wedge \tau_0^-)} v_y(X_{t \wedge \tau_0^-}) : t \geq 0\}$ is a $\mathbb{P}_{x-y}$-martingale. Hence, $\{e^{-(t \wedge \tau_y^-)} v_y(X_{t \wedge \tau_y^-}) : t \geq 0\}$ is a $\mathbb{P}_x$-martingale.

*Supermartingale property* (ii). On the event $\{t < \mathbf{e}_r\}$, the identity $\underline{X}_{\mathbf{e}_r} = \underline{X}_t \wedge (X_t + I)$ holds and, conditionally on $\mathcal{F}_t$, $I$ has the same distribution as $\underline{X}_{\mathbf{e}_r}$. In particular, it follows that, on $\{t < \mathbf{e}_r\}$, $\underline{X}_{\mathbf{e}_r} \leq X_t + I$. If

$$e^y \leq K\mathbb{E}(e^{\underline{X}_{\mathbf{e}_r}}), \tag{14}$$

then, for $x \in \mathbb{R}$,

$$v_y(x) \geq \frac{\mathbb{E}[\mathbf{1}_{(t < \mathbf{e}_r)}\mathbb{E}[(K\mathbb{E}[e^{x+\underline{X}_{\mathbf{e}_r}}] - e^{X_t + I})\mathbf{1}_{(-(X_t+I) > x-y)}|\mathcal{F}_t]]}{\mathbb{E}[e^{\underline{X}_{\mathbf{e}_r}}]}$$
$$\geq \mathbb{E}[e^{-rt} v_y(x + X_t)]$$
$$= \mathbb{E}_x[e^{-rt} v_y(X_t)].$$

Using the stationary independent increments of $X$, the latter inequality is sufficient to deduce that $\{e^{-rt} v_y(X_t) : t \geq 0\}$ is a $\mathbb{P}_x$-supermartingale.

*Continuous pasting* (iii). It is clear from (13) that, under (14), $v_y(x) \geq 0$ for all $x \in \mathbb{R}$. It is also a straightforward manipulation to show that

$$(15) \qquad v_y(x) = (K - e^x) + \frac{\mathbb{E}[(e^{x+\underline{X}_{\mathbf{e}_r}} - K\mathbb{E}(e^{\underline{X}_{\mathbf{e}_r}}))\mathbf{1}_{(-\underline{X}_{\mathbf{e}_r} \leq x - y)}]}{\mathbb{E}[e^{\underline{X}_{\mathbf{e}_r}}]}.$$

When $x < y$ and (14) holds, we see that $v_y(x) = (K - e^x) = (K - e^x)^+$, as required. When $x = y$, in order to satisfy condition (i), we require that

$$\frac{e^y - K\mathbb{E}[e^{\underline{X}_{\mathbf{e}_r}}]}{\mathbb{E}[e^{\underline{X}_{\mathbf{e}_r}}]}\mathbb{P}(-\underline{X}_{\mathbf{e}_r} = 0) = 0.$$

If $\mathbb{P}(-\underline{X}_{\mathbf{e}_r} = 0) = 0$, that is to say, if 0 is regular for $(-\infty, 0)$, then continuous pasting always holds. If, on the other hand, $\mathbb{P}(-\underline{X}_{\mathbf{e}_r} = 0) > 0$, that is to say 0 is irregular for $(-\infty, 0)$, then we are obliged to choose $y = x^*$; in other words, $e^y = K\mathbb{E}[e^{\underline{X}_{\mathbf{e}_r}}]$.

*Lower bound* (iv). From (15), we see in the case that 0 is irregular for $(-\infty, 0)$, with the established choice $y = x^*$ for continuity, that automatically $v_y(x) \geq (K - e^x)$. Hence, together with $v_y(x) \geq 0$, the lower bound (ii) is satisfied.

When 0 is regular for $(-\infty, 0)$, we know that if (14) holds, then $v_y(x) \geq 0$. On the other hand, considering (15), we see that a sufficient condition that $v_y(x) \geq (K - e^x)$ is that

$$(16) \qquad e^y \geq K\mathbb{E}[e^{\underline{X}_{\mathbf{e}_r}}].$$

Since $v_y(x) \geq (K - e^x)^+$ if and only if $v_y(x) \geq 0$ and $v_y(x) \geq (K - e^x)$, then it becomes clear from (14) and (16) that, to respect the inequality (ii), we may choose $y = x^*$. There can be no other choice of $y$ with this property, for then we would be able to produce two distinct solutions to the optimal stopping problem. $\square$

REMARK 9. In the above proof, by assuming from the outset that the optimal threshold is equal to $x^*$, the proof takes the form of a direct verification of Mordecki's solution. However, in that case, one does not get a feel for the role played by regularity in the solution.

PROOF OF THEOREM 6. Since $v(x) = K - e^x$ for all $x \leq x^*$, and, hence, $v'(x^*-) = -e^{x^*}$, we are required to show that $v'(x^*+) = -e^{x^*}$ for smooth pasting. Starting from (13) and recalling that $e^{x^*} = K\mathbb{E}[e^{\underline{X}_{\mathbf{e}_r}}]$, we have

$$v(x) = -K\mathbb{E}[(e^{x-x^*+\underline{X}_{\mathbf{e}_r}} - 1)\mathbf{1}_{(-\underline{X}_{\mathbf{e}_r} > x - x^*)}]$$
$$= -K(e^{x-x^*} - 1)\mathbb{E}[e^{\underline{X}_{\mathbf{e}_r}}\mathbf{1}_{(-\underline{X}_{\mathbf{e}_r} > x - x^*)}]$$
$$\quad - K\mathbb{E}[(e^{\underline{X}_{\mathbf{e}_r}} - 1)\mathbf{1}_{(-\underline{X}_{\mathbf{e}_r} > x - x^*)}].$$

From the last equality, we may then write

$$\frac{v(x) - (K - e^{x^*})}{x - x^*} = \frac{v(x) + K(\mathbb{E}[e^{\underline{X}_{\mathbf{e}_r}}] - 1)}{x - x^*}$$

$$= -K\frac{(e^{x-x^*} - 1)}{x - x^*}\mathbb{E}[e^{\underline{X}_{\mathbf{e}_r}}\mathbf{1}_{(-\underline{X}_{\mathbf{e}_r} > x - x^*)}]$$

$$+ K\frac{\mathbb{E}[(e^{\underline{X}_{\mathbf{e}_r}} - 1)\mathbf{1}_{(-\underline{X}_{\mathbf{e}_r} \leq x - x^*)}]}{x - x^*}.$$

To simplify notation, let us call $A_x$ and $B_x$ the last two terms, respectively. It is clear that

(17) $$\lim_{x \downarrow x^*} A_x = -K\mathbb{E}[e^{\underline{X}_{\mathbf{e}_r}}\mathbf{1}_{(-\underline{X}_{\mathbf{e}_r} > 0)}].$$

On the other hand, we have that

$$B_x = K\frac{\mathbb{E}[(e^{\underline{X}_{\mathbf{e}_r}} - 1)\mathbf{1}_{(0 < -\underline{X}_{\mathbf{e}_r} \leq x - x^*)}]}{x - x^*}$$

$$= K\int_{0+}^{x-x^*} \frac{e^{-z} - 1}{x - x^*}\mathbb{P}(-\underline{X}_{\mathbf{e}_r} \in dz)$$

$$= K\frac{e^{x^* - x} - 1}{x - x^*}\mathbb{P}(0 < -\underline{X}_{\mathbf{e}_r} \leq x - x^*)$$

$$+ \frac{K}{x - x^*}\int_{0+}^{x-x^*} e^{-z}\mathbb{P}(0 < -\underline{X}_{\mathbf{e}_r} \leq z)\,dz,$$

where in the first equality we have removed the possible atom at zero from the expectation by noting that $\exp\{\underline{X}_{\mathbf{e}_r}\} - 1 = 0$ on $\underline{X}_{\mathbf{e}_r} = 0$. This leads to $\lim_{x \downarrow x^*} B_x = 0$. Using the expression for $e^{x^*}$, we see that $v'(x^*+) = -e^{x^*} + K\mathbb{P}(-\underline{X}_{\mathbf{e}_r} = 0)$, which equals $-e^{x^*}$ if and only if $\mathbb{P}(-\underline{X}_{\mathbf{e}_r} = 0) = 0$; in other words, if and only if 0 is regular for $(-\infty, 0)$. □

REMARK 10. The intuition behind the role of regularity in the solution to the optimal stopping problem is as follows.

For each $y$, the function $v_y$ provides us with a martingale up to stopping at $\tau_y^-$ and, otherwise, a supermartingale providing $y$ is no greater than $K\mathbb{E}[e^{\underline{X}_{\mathbf{e}_r}}]$. We then appeal to the classical interpretation of the value function of the optimal stopping problem as the least superharmonic majorant of the gain.

If there is irregularity, then choosing a candidate $y$ for the optimal level which is strictly smaller than $K\mathbb{E}[e^{\underline{X}_{\mathbf{e}_r}}]$ allows for the function $v_y$ to be lower bounded by the gain, but introduces a discontinuity at $y$. The closer to $K\mathbb{E}[e^{\underline{X}_{\mathbf{e}_r}}]$ we bring $y$ from above, the closer we bring $v_y$ pointwise to the gain function. Choosing the candidate level $y$ below $K\mathbb{E}[e^{\underline{X}_{\mathbf{e}_r}}]$ cannot guarantee that $v_y$ will be bounded below by the gain and, hence, $\tau_y^-$ is not an admissible strategy. Here we experience a principle of continuous pasting in order to fulfill the traditional role of the gain function as being the smallest superharmonic majorant of the gain.

If there is regularity, then all curves $v_y$ are continuous. For $y$ less than $K\mathbb{E}[e^{\underline{X}_{\mathbf{e}_r}}]$ there is, in general, a discontinuity in the first derivative of $v_y$ at the point $y$. As we move $y$ closer to $K\mathbb{E}[e^{\underline{X}_{\mathbf{e}_r}}]$ from above, the discrepancy in the first derivative at $y$ tends to zero and the function $v_y$ moves pointwise closer to the gain function. The functions $v_y$ for which $y$ is strictly greater than $K\mathbb{E}[e^{\underline{X}_{\mathbf{e}_r}}]$ do not bound the gain from above and, hence, again offer inadmissible strategies. We therefore experience a principle of smooth pasting this time in order to fulfill the traditional role of the gain function as being the smallest superharmonic majorant of the gain.

This explaination ties up to the remarks of Boyarchenko and Levendorskiĭ [12], who talk about the optimal value function being "more smooth" than other potential candidates.

REMARK 11. It is also interesting to note that, for the case that, 0 is irregular for $(-\infty, 0)$, the conditions (i)–(iii) suffice to solve the optimal stopping problem. In this sense, the system (i)–(iv) is "closed" by the principle of continuous pasting. In the regular case the extra condition (iv) is genuinely required in order to "close" the system (i)–(iv) and establish the unique solution to (1). Hence, in the regular case, smooth pasting is implicitly required for closure of the system (i)–(iv).

**7. Conclusions and conjectures.** We have emphasized a simple fluctuation identity which binds together a number of other known identities for general and special choices of Lévy processes. Our identity helps to explain the nature of known solutions to the American put optimal stopping problem, giving the opportunity to explore the pasting principle in more detail. Further, we have extracted the precise criterion for the phenomenon of smooth pasting versus continuous pasting from Mordecki's solution to (1). We conjecture that this criterion remains unchanged when moving to the same optimal stopping problem, but with finite horizon. In this case, simulations of Matache, Nitsche and Schwab [26], Hirsa and Madan [23], Almendral and Oosterlee [2] and Almendral [1] show that the optimal strategy consists of stopping when first crossing below a boundary which is nonincreasing in time. Therefore, since the boundary has only bounded variation, and assuming it is also right continuous, the regularity of the boundary for the stopping domain below it should correspond to the regularity of 0 for $(-\infty, 0)$ and the same necessary and sufficient conditions should dictate the pasting condition.

In general, it would seem that one can use regularity as a reasonable rule of thumb for predicting when smooth pasting occurs. Roughly speaking, suppose one has an optimal stopping problem for a Markov process $X$ with expectation operators $\{\mathbf{E}_x\}$ taking the form

$$u(x) = \sup_{\tau \in \mathcal{T}} \mathbf{E}_x[G(X_\tau)],$$

where the family of stopping times $\mathcal{T}$ and smooth function $G$ are such that $u$ is well defined. Assume that the optimal strategy corresponds to the partitioning of the domain of $X$ into two regions, say, $\mathcal{C}$ and $\mathcal{S} = \mathcal{C}^c$, and stopping as soon as $X$ enters $\mathcal{S}$. (Note that so far all of these assumptions conform to most solutions which are known in explicit or semi-explicit form.) It would then be reasonable to work with the assumption that there will be smooth pasting at $x \in \partial \mathcal{S}$ if and only if $x$ is regular for $\mathcal{S}$; otherwise continuous pasting holds. Note that the exercise boundary may consist of disjoint regions which may exhibit different pasting principles. At the very least, in addition to the articles discussed in this paper concerning the American put optimal stopping problem, pasting principles based on path regularity would also seem to apply to the optimal stopping problems in [32, 33] and [20] in the context of a sequential testing and Poisson disorder problems and again in [4] and [6] for the case of pricing perpetual Russian options under Lévy models.

**Acknowledgments.** This work was carried out during visits of the first author to Utrecht University and the second author to ETH Zürich and MaPhySto at Aarhus University. Both would like to thank the host institutions of their counterpart for their hospitality and support and, in addition, A. E. Kyprianou is likewise grateful to MaPhySto. Thanks are due to Ana–Maria Matache, Christoph Schwab, Ron Doney, Sergei Levendorskiĭ and Freddy Delbaen for a number of helpful discussions. Special thanks go to Goran Peskir who, through a number of insightful discussions, has inspired this work and greatly helped us to understand the *principle* of continuous pasting.

DEPARTMENT OF STATISTICS
THE UNIVERSITY OF WARWICK
COVENTRY CV4 7AL
UNITED KINGDOM
E-MAIL: L.Alili@warwick.ac.uk

UTRECHT UNIVERSITY
P.O. BOX 80.010
3500 TA, UTRECHT
THE NETHERLANDS
E-MAIL: kyprianou@math.uu.nl